\documentclass[reqno]{article}
\usepackage{amsmath}
\usepackage{amssymb}
\def\Re{\mathop{\mathrm{Re}}}

\begin{document}

\title{\bf An order result for the exponential divisor function}
\author{{\sc L\'aszl\'o T\'oth} \thanks{Research supported by the fund of
Applied Number Theory Research Group of the Hungarian Academy of
Sciences.} \\
University of P\'ecs, Institute of Mathematics and Informatics,\\
Ifj\'us\'ag u. 6, 7624 P\'ecs, Hungary,\\
ltoth@ttk.pte.hu}
\date{ }
\maketitle

\centerline{Publ. Math. Debrecen, {\bf 71} (2007), no. 1-2, 165-171}
\vskip4mm



{\it Abstract} : The integer $d=\prod_{i=1}^s p_i^{b_i}$ is called
an exponential divisor of $n=\prod_{i=1}^s p_i^{a_i}>1$ if $b_i \mid
a_i$ for every $i\in \{1,2,...,s\}$. Let $\tau^{(e)}(n)$ denote the
number of exponential divisors of $n$, where $\tau^{(e)}(1)=1$ by
convention. The aim of the present paper is to establish an
asymptotic formula with remainder term for the $r$-th power of the
function $\tau^{(e)}$, where $r\ge 1$ is an integer. This improves
an earlier result of {\sc M. V. Subbarao} \cite{Su72}.

\vskip7mm {\bf 1. Introduction} \vskip3mm

Let $n>1$ be an integer of canonical form $n=\prod_{i=1}^s
p_i^{a_i}$. The integer $d$ is called an {\sl exponential divisor}
of $n$ if $d=\prod_{i=1}^s p_i^{b_i}$, where $b_i \mid a_i$ for
every $i\in \{1,2,...,s\}$, notation: $d\mid_e n$. By convention
$1\mid_e 1$.

Let $\tau^{(e)}(n)$ denote the number of exponential divisors of
$n$. The function $\tau^{(e)}$ is called the {\sl exponential
divisor function}. {\sc J. Wu} \cite{Wu95} showed, improving an
earlier result of {\sc M. V. Subbarao} \cite{Su72}, that
$$
\sum_{n\le x} \tau^{(e)}(n)=Ax +B x^{1/2} + O(x^{2/9}\log x),
\leqno(1)
$$
where
$$
A: =\prod_p \left(1 + \sum_{a=2}^{\infty}
\frac{\tau(a)-\tau(a-1)}{p^a} \right), \, B: =\prod_p \left(1 +
\sum_{a=5}^{\infty}
\frac{\tau(a)-\tau(a-1)-\tau(a-2)+\tau(a-3)}{p^{a/2}} \right),
$$
$\tau$ denoting the usual divisor function. The $O$-term can further
be improved.

Other properties of the function $\tau^{(e)}$, compared with those
of the divisor function $\tau$ were investigated in papers
\cite{KaSu2003}, \cite{KoIv97}, \cite{SmWu97}, \cite{Su72}.

{\sc M. V. Subbarao} \cite{Su72} remarked that for every positive
integer $r$,
$$
\sum_{n\le x} (\tau^{(e)}(n))^r \sim A_rx, \leqno(2)
$$
where
$$
A_r: =\prod_p \left(1 + \sum_{a=2}^{\infty}
\frac{(\tau(a))^r-(\tau(a-1))^r}{p^a} \right). \leqno(3)
$$

It is the aim of the present paper to establish the following more
precise asymptotic formula for the $r$-th power of the function
$\tau^{(e)}$, where $r\ge 1$ is an integer:
$$
\sum_{n\le x} (\tau^{(e)}(n))^r = A_r x + x^{1/2} P_{2^r-2}(\log x)
+ O(x^{u_r+\varepsilon}), \leqno(4)
$$
for every $\varepsilon >0$, where $A_r$ is given by (3), $P_{2^r-2}$
is a polynomial of degree $2^r-2$ and $u_r:
=\frac{2^{r+1}-1}{2^{r+2}+1}$.

Note that a similar formula is known for the divisor function $\tau$, namely for
any integer $r\ge 2$,
$$
\sum_{n\le x} (\tau(n))^r =  x Q_{2^r-1}(\log x) +
O(x^{v_r+\varepsilon}), \leqno(5)
$$
valid for every $\varepsilon >0$, where $v_r: =\frac{2^r-1}{2^r+2}$
and $Q_{2^r-1}$ is a polynomial of degree $2^r-1$, this goes back to
the work of {\sc S. Ramanujan}, cf. \cite{Wi22}.

Formula (4) is a direct consequence of a simple general result,
given in Section 2 as Theorem, regarding certain multiplicative
functions $f$ such that $f(n)$ depends only on the $\ell$-full
kernel of $n$, where $\ell \ge 2$ is a fixed integer.

We also consider a generalization of the exponential divisor
function, see Section 4.

Let $\phi^{(e)}(n)$ denote the number of divisors $d$ of $n$ such
that $d$ and $n$ have no common exponential divisors. The function
$\phi^{(e)}$ is multiplicative and for every prime power $p^a$
($a\ge 1$), $\phi^{(e)}(p^a)= \phi(a)$, where $\phi$ is the Euler
function.

As another consequence of our Theorem we obtain for every integer
$r\ge 1$ that
$$
\sum_{n\le x} (\phi^{(e)}(n))^r = B_rx + x^{1/3} R_{2^r-2}(\log x) +
O(x^{t_r+\varepsilon}), \leqno(6)
$$
for every $\varepsilon >0$, where $t_r: =\frac{2^{r+1}-1}{3\cdot
2^{r+1}}$, $R_{2^r-2}$ is a polynomial of degree $2^r-2$ and
$$
B_r: =\prod_p \left(1 + \sum_{a=3}^{\infty}
\frac{(\phi(a))^r-(\phi(a-1))^r}{p^a} \right). \leqno(7)
$$

In the case $r=1$ formula (6) was proved in \cite{To2004} with a
better error term.

Our error terms depend on estimates for
$$
D(1,\underbrace{\ell,\ell,...,\ell}_{k-1};x): = \sum_{ab_1^{\ell}b_2^{\ell}\cdot ...
\cdot b_{k-1}^{\ell}\le x} 1,
$$
where $k, \ell \ge 2$ are fixed and $a, b_1, b_2,\cdots ,b_{k-1}\ge 1$ are integers.
\vskip3mm {\bf 2. A general result} \vskip3mm

We prove the following general result.

\vskip2mm {\bf Theorem.} Let $f$ be a complex valued multiplicative
arithmetic function such that

a) $f(p)=f(p^2)=\cdots =f(p^{\ell-1})=1$,
$f(p^{\ell})=f(p^{\ell+1})=k$ for every prime $p$, where $\ell, k\ge
2$ are fixed integers and

b) there exist constants $C, m>0$ such that $|f(p^a)|\le Ca^m$ for
every prime $p$ and every $a\ge \ell +2$.

Then for $s\in \mathbb{C}$

i) $$ F(s):= \sum_{n=1}^{\infty} \frac{f(n)}{n^s} =
\zeta(s)\zeta^{k-1}(\ell s) V(s), \qquad \Re s > 1,
$$
where the Dirichlet series $V(s):=\sum_{n=1}^{\infty} \frac{v(n)}{n^s}$ is absolutely
convergent for $\Re s > \frac1{\ell+2}$,

ii) $$\sum_{n\le x} f(n)= C_fx + x^{1/\ell} P_{f,k-2}(\log x)
+O(x^{u_{k,\ell}+\varepsilon}),
$$
for every $\varepsilon >0$, where $P_{f,k-2}$ is a polynomial of
degree $k-2$, $u_{k,\ell}: = \frac{2k-1}{3+(2k-1)\ell}$ and
$$
C_f: =\prod_p \left(1 + \sum_{a=\ell}^{\infty}
\frac{f(p^a)-f(p^{a-1})}{p^a} \right).
$$

iii) The error term can be improved for certain values of $k$ and
$\ell$. For example in the case $k=3$, $\ell=2$ it is
$O(x^{8/25}\log^3 x)$.

\vskip3mm
{\bf 3. Proofs}
\vskip3mm

The proof of the Theorem is based on the following Lemma.
For an integer $\ell \ge 1$ let $\mu_{\ell}(n)=\mu(m)$ or $0$, according as
$n=m^{\ell}$ or not, where $\mu$ is the M\"obius function. Note that function
$\mu_{\ell}$ is multiplicative and for any prime power $p^a$ ($a\ge 1$),
$$
\mu_{\ell}(p^a)= \begin{cases} -1, & \text{ if } a=\ell, \\
                          0, & \text{ otherwise.} \end{cases} \leqno(8)
$$
Furthermore, for an integer $h\ge 1$ let the function
$\mu_{\ell}^{(h)}$ be defined in terms of the Dirichlet convolution
by
$$
\mu_{\ell}^{(h)}=\underbrace{\mu_{\ell} * \mu_{\ell}*\cdots *\mu_{\ell}}_h.
$$
The function $\mu_{\ell}^{(h)}$ is also multiplicative.

\vskip2mm
{\bf Lemma.} For any integers $h, \ell \ge 1$ and any prime power $p^a$ ($a\ge 1$),
$$
\mu_{\ell}^{(h)}(p^a)= \begin{cases} (-1)^j {h \choose j}, & \text{
if } a=j\ell, \quad 1\le j\le h, \\
                                 0, & \text{ otherwise. } \end{cases} \leqno(9)
$$
\vskip2mm

\vskip2mm {\bf Proof of the Lemma.} By induction on $h$. For $h=1$
this follows from (8). We suppose that formula (9) is valid for $h$
and prove it for $h+1$. Using the relation
$\mu_{\ell}^{(h+1)}=\mu_{\ell}^{(h)}* \mu_{\ell}$ and (8) we obtain
for $a<\ell$,
$$
\mu_{\ell}^{h+1}(p^a)=\mu_{\ell}^{(h)}(p^a)=0
$$
and for $a\ge \ell$,
$$
\mu_{\ell}^{(h+1)}(p^a)=\mu_{\ell}^{(h)}(p^a)- \mu_{\ell}^{h}(p^{a-\ell})
$$
$$
= \begin{cases} \mu_{\ell}^{(h)}(p^{\ell})- 1 = - {h \choose 1} -1 =
-{h+1 \choose 1}, & \text{ if } a=\ell, \\
(-1)^j {h \choose j} -(-1)^{j-1} {h \choose j-1}= (-1)^j {h+1
\choose j}, & \text{ if } a=\ell j, \quad 2\le j \le h, \\
-\mu_{\ell}^{(h)}(p^{h\ell}) = -(-1)^h {h \choose h} = (-1)^{h+1}
{h+1 \choose h+1}, & \text{ if } a=(h+1)\ell, \\ 0, & \text{
otherwise, } \end{cases}
$$
which proves the Lemma. \vskip2mm

\vskip2mm {\bf Proof of the Theorem.} i) We can formally obtain the
desired expression by taking $v=f*\mu * \mu_{\ell}^{(k-1)}$. Here
$v$ is multiplicative and easy computations show that $v(p^a)=0$ for
any $1\le a\le \ell +1$ and for $a\ge \ell +2$,
$$
v(p^a)= \sum_{j\ge 0} (-1)^j {k-1 \choose j}
\left(f(p^{a-j\ell})-f(p^{a-j\ell-1})\right),
$$
where, according to the Lemma, the number of nonzero terms is at
most $k$.

Let $M_k=\max_{0\le j\le k-1} {k-1 \choose j}$. We obtain that for
every prime $p$ and every $a\ge \ell +2$,
$$
|v(p^a)|\le 2k M_k C a^m.
$$
For every $\varepsilon >0$, $a^m\le 2^{a\varepsilon}$ for sufficiently large $a$,
$a\ge a_0$ say, where $a_0\ge \ell +2$. For $\Re s > 1/(\ell+2)$ choose $\varepsilon>
0$ such that $\Re s - \varepsilon > 1/(\ell+2)$. Then
$$
\sum_p \sum_{a\ge a_0} \frac{|v(p^a)|}{p^{as}}\le 2k M_k C \sum_p \sum_{a\ge a_0}
\frac{2^{a\varepsilon}}{p^{as}} \le  2k M_k C \sum_p \sum_{a\ge a_0}
\frac1{p^{a(s-\varepsilon)}}=
$$
$$
= 2k M_k C \sum_p
\frac1{p^{a_0(s-\varepsilon)}}\left(1-\frac1{p^{s-\varepsilon}}\right)^{-1} \le 2k
M_k C \left(1-\frac1{2^{1/(\ell+2)}}\right)^{-1} \sum_p
\frac1{p^{a_0(s-\varepsilon)}},
$$
and obtain that $V(s)$ is absolutely convergent for $Re\; s
> 1/(\ell+2)$.

Note that $v(p^{\ell+2})=f(p^{\ell+2})-k$ for every $\ell \ge 3, k\ge 2$ and for
$\ell=2, k\ge 2$ it is $v(p^4)=f(p^4)-{k+1 \choose 2}$.

ii) Consider the $k$-dimensional generalized divisor function
$$
d(1,\underbrace{\ell,\ell,...,\ell}_{k-1};n)=
\sum_{ab_1^{\ell}b_2^{\ell}\cdots b_{k-1}^{\ell}=n} 1.
$$
According to i),
$$
f(n)=\sum_{ab=n} d(1,\underbrace{\ell,\ell,...,\ell}_{k-1};a) v(b).
$$
One has, see \cite{Kr88}, Ch. 6,
$$
\sum_{n\le x} d(1,\underbrace{\ell,\ell,...,\ell}_{k-1};n)=
\leqno(10)
$$ $$
=K_1 x + x^{1/\ell} \left(K_2 \log^{k-2} x + K_3 \log^{k-3} x +
\cdots + K_{k-1} \log x + K_k\right) + O(x^{u_{k,\ell}
+\varepsilon}),
$$
for every $\varepsilon >0$, where $u_{k,\ell}=
\frac{2k-1}{3+(2k-1)\ell}$ (see \cite{Kr88}, Theorem 6.10),
$K_1,K_2,...,K_{k-1},K_k$ are absolute constants depending on $k$
and $\ell$ and $K_1=\zeta^{k-1}(\ell)$. For example for $k=2$ one
has $K_2=\zeta(\frac1{\ell})$, and for $k=3$:
$K_2=\frac1{\ell}\zeta(\frac1{\ell})$,
$K_3=(2\gamma-1)\zeta(\frac1{\ell}) +\frac1{\ell}
\zeta'(\frac1{\ell})$, where $\gamma$ is Euler's constant.

We obtain
$$
\sum_{n\le x} f(n) =\sum_{ab\le x}
d(1,\underbrace{\ell,\ell,...,\ell}_{k-1};a) v(b) =\sum_{b\le x}
v(b) \sum_{a\le x/b} d(1,\underbrace{\ell,\ell,...,\ell}_{k-1};a) =
$$
$$ = \sum_{b\le x} v(b) \left( K_1 (x/b) + (x/b)^{1/\ell} \left( K_2
\log^{k-2} (x/b) + K_3 \log^{k-3} (x/b) + \cdots \right. \right.
$$ $$\left. \left. + K_{k-1} \log (x/b) + K_k \right) + O( (x/b)^{u +\varepsilon})
\right),
$$
and obtain the desired result by partial summation and by noting
that $u_{k,\ell} > 1/(\ell+2)$.

iii) For $k=3$, $\ell=2$ the errror term of (10) is
$O(x^{8/25}\log^3 x)$, cf. \cite{Kr88}, Theorem 6.4.

\vskip3mm

{\bf 4. Applications.}  1. In case $f(n)= (\tau^{(e)}(n))^r$, where
$r\ge 1$ is an integer, we obtain formula (4) applying the Theorem
for $\ell=2$, $k=2^r$.

\vskip2mm

2. For $k\ge 2$ consider the multiplicative function $f(n)=\tau_k^{(e)}(n)$, where
for every prime power $p^a$ ($a\ge 1$), $\tau_k^{(e)}(p^a) :=\tau_k(a)$ representing
the number of ordered $k$-tuples of positive integers $(x_1,...,x_k)$ such that
$a=x_1\cdot ... \cdot x_k$. Here $\tau_k(p^b)= {b+k-1 \choose k-1}$ for every prime
power $p^b$ ($b\ge 1$). In case $k=2$, $\tau_2^{(e)}(n)= \tau^{(e)}(n)$.

Taking $\ell=2$ and $k:=k$ we obtain that $v(p^4)=\tau_k(4)-k(k+1)/2=0$ and $V(s)$ is
absolutely convergent for $\Re s > \frac1{5}$ (and not only for $\Re s>\frac1{4}$
given by the Theorem),
$$
\sum_{n\le x} \tau_k^{(e)}(n)= C_k x + x^{1/2}
S_{k-2}(\log x) + O(x^{w_k+\varepsilon}), \leqno(11)
$$
for every $\varepsilon >0$, where $S_{k-2}$ is a polynomial of
degree $k-2$, $w_k: =\frac{2k-1}{4k+1}$ and
$$
C_k=\prod_p \left(1 + \sum_{a=2}^{\infty}
\frac{\tau_k(a)-\tau_k(a-1)}{p^a} \right).
$$
For $k=3$ the error term of (11) can be improved into
$O(x^{8/25}\log^3 x)$.

A similar formula can be obtained for $\sum_{n\le x}
(\tau_k^{(e)}(n))^r$. \vskip2mm

3. For the function $\phi^{(e)}(n)$ defined in the Introduction we
obtain formula (6) by choosing $\ell=3$, $k=2^r$.

\vskip4mm


\vskip3mm


\begin{thebibliography}{99}

\bibitem{KaSu2003} {\sc I. K\'atai} and {\sc M. V. Subbarao}, On the
distribution of exponential divisors, {\it Annales Univ. Sci. Budapest.,
Sect. Comp.}, {\bf 22} (2003), 161-180.

\bibitem{KoIv97} {\sc J. -M. de Koninck} and {\sc A. Ivi\'c}, An asymptotic formula
for reciprocals of logarithms of certain multiplicative functions
{\it Canad. Math. Bull.}, {\bf 21} (1978), 409-413.

\bibitem{Kr88} {\sc E. Kr\"atzel}, {\it Lattice points}, Kluwer, Dordrecht-Boston-London, 1988.

\bibitem{SmWu97} {\sc A. Smati} and {\sc J. Wu}, On the exponential divisor function,
{\it Publ. Inst. Math. (Beograd) (N. S.)}, {\bf 61} (1997), 21-32.

\bibitem{Su72} {\sc M. V. Subbarao}, On some arithmetic convolutions, in {\it The Theory of Arithmetic Functions}, Lecture Notes in Mathematics No. {\bf 251},
247-271, Springer, 1972.

\bibitem{To2004} {\sc L. T\'oth}, On certain arithmetic functions involving
exponential divisors, {\it Annales Univ. Sci. Budapest., Sect.
Comp.}, {\bf 24} (2004), 285-294.

\bibitem{Wu95} {\sc J. Wu}, Probl\`eme de diviseurs exponentiels et entiers
exponentiellement sans facteur carr\'e, {\it J. Th\'eor. Nombres
Bordeaux}, {\bf 7} (1995), 133-141.

\bibitem{Wi22} {\sc B. M. Wilson}, Proofs of some formulae enunciated by Ramanujan,
{\it Proc. London Math. Soc. (2)}, {\bf 21} (1922), 235-255.

\end{thebibliography}
\end{document}